\documentclass[3p, nopreprintline]{elsarticle}

\usepackage{hyperref}
\usepackage{units}
\usepackage{gensymb}
\usepackage{pgfplots}
\usepackage{amsmath}
\usepackage{txfonts}
\usepackage{bm}

\usepackage{caption}
\usepackage{subcaption}

\usepackage{todonotes}
\usepackage{layouts}
\usepackage{siunitx}

\usepackage[section]{placeins}

\usetikzlibrary{fadings}
\usetikzlibrary{patterns}
\usetikzlibrary{shadows.blur}
\usetikzlibrary{shapes}

\journal{Finite Elements in Analysis and Design}

\DeclareUnicodeCharacter{2212}{-}
\begin{document}

\begin{frontmatter}

\title{Numerical Design of Distributive Mixing Elements}

%% Group authors per affiliation:
%\author{Elsevier\fnref{myfootnote}}
%\address{Radarweg 29, Amsterdam}
%\fntext[myfootnote]{Since 1880.}

%% or include affiliations in footnotes:
\author[catsaddress]{Sebastian Hube\corref{mycorrespondingauthor}}
\cortext[mycorrespondingauthor]{Corresponding author}
\ead{eusterholz@cats.rwth-aachen.de}

\author[catsaddress]{Marek Behr}
\ead{behr@cats.rwth-aachen.de}

\author[catsaddress,ilsbaddress]{Stefanie Elgeti}
\ead{elgeti@cats.rwth-aachen.de}

\author[ikvaddress]{Malte Schön}
\ead{malte.schoen@ikv.rwth-aachen.de}

\author[ikvaddress]{Jana Sasse}
\ead{jana.sasse@ikv.rwth-aachen.de}

\author[ikvaddress]{Christian Hopmann}
\ead{office@ikv.rwth-aachen.de}

\address[catsaddress]{Chair for Computational Analysis of Technical Systems (CATS), RWTH Aachen University, Schinkelstr. 2, 52062 Aachen, Germany}
\address[ikvaddress]{Institute for Plastics Processing (IKV) at the RWTH Aachen University, Seffenter Weg 201, 52074 Aachen, Germany}
\address[ilsbaddress]{Institute of Lightweight Design and Structural Biomechanics (E317), TU Wien, Gumpendorfer Str. 7, A-1060 Vienna, Austria}

\begin{abstract}
This paper presents a novel shape-optimization technique for the design of mixing elements in single-screw extruders.
Extruders enable the continuous production of constant-cross-section profiles.
Equipped with one or several screw-shaped rotors, the extruder transports solid polymer particles towards the outlet.
Due to shear heating, melting is induced and a melt stream is created, which can be further processed.
While many variants of multi-screw extruders exist, a significant share of all extrusion machines is made up of single-screw extruders due to their comparatively low operating costs and complexity.
While the reduced complexity yields economic benefits, single-screw extruders' mixing capabilities, i.e., their ability to produce a melt with a  homogeneous material and temperature distribution, suffer compared to multi-screw extruders.
To compensate for this shortcoming, so-called mixing elements are added to the screw to enhance dynamic mixing by recurring flow reorientations.
In view of the largely unintuitive flow characteristics of polymer melts, we present an optimization framework that allows designing these mixing elements numerically based on finite-element simulations of the melt flow.
To reduce the computational demand required by shape optimization of a complete mixing section, we only focus on the shape optimization of a single mixing element.
This paper presents advances in three aspects of numerical design:
(1) A combination of free-form deformation and surface splines is presented, allowing to parameterize the mixing element's shape by very few variables.
(2) The combination of this concept with a linear-elasticity-based mesh update method to deform the computational domain without the need for remeshing is demonstrated.
(3) A simple yet robust and sensitive new objective formulation to assess distributive mixing in laminar flows based on a measure for the interfacial area is proposed for the optimization.

\end{abstract}


\end{frontmatter}

\section{Introduction}
In plastics processing, it is very common to generate plastics melt in screw-based units. 
%Within the many ways of plastics processing, extrusion is one of the most prominent techniques.
One example are single-screw extruders, which are used to provide constant and homogeneous melt streams within profile extrusion.
% streams, which in profile extrusion are further processed by conveying the melt through a die.
As another example, we refer to the injecting molding, and its screw-based plasticating units.
%To obtain such homogeneous melt streams, the extruder and the plasticating unit are  continuously fed with solid polymers.
Starting point are always solid polymer pellets, which are continuously inserted into the screw unit, whose purpose is to melt, mix, and pump the polymer material \cite{Erwin78_TheoryOfMixingSectionsInSingleScrewExtruders}.
%However, in injection modling, the melt is pushed into a mold by a translational movement of the screw.
For screw units, two design categories are popular: single-screw and twin-screw versions. 
It can be observed that single-screw extruders are more commonly used. 
This is mainly due to their easy and economical installation and operation.
However, they also suffer from the drawback that the provided melt stream is inferior to twin-screw extruders in terms of homogeneity of both physical properties of the base melt or blends and distribution of additives.
This often has a negative impact on the quality of the final plastics product.
This situation is amplified by the fact that many additives can only be mixed into the fully molten material; and therefore are added towards the end of the extruder screw.
In order to still achieve a homogeneous distribution of the additive within the melt, special screw elements that are responsible for mixing are included in this area. \par
Depending on the desired effects, we frequently distinguish between two types of mixing: dispersive and distributive mixing.
While dispersive mixing describes breaking down agglomerates reducing particle clusters, distributive mixing describes homogenizing the particle distribution within the melt.
In  this work, we focus on enhancing distributive mixing by numerical design.\par
While achieving good mixing in turbulent flows is a comparably well-understood problem, obtaining similar effects in laminar flows -- exactly the flow type predominant in extruders due to the high viscosity of polymer melts -- is more challenging.
In particular, within extruders, diffusion and turbulence are negligible, and shear remains the only dominant mixing mechanism \cite{Mohr57}.
Consequently, for the design of extruders, the theory of laminar mixing is particularly relevant \cite{Mohr57, Spencer51, Erwin78_TheoryOfLaminarMixing}.
Based on these findings, we can derive the design principles for extruder screws in general, and in particular for the mixing sections \cite{Erwin78_TheoryOfMixingSectionsInSingleScrewExtruders, Erwin78_AnUpperBoundOnThePerformanceOfPLaneStrainMixers}.
In these mixing sections, obstacles -- often cylindrically or rhomboidally shaped -- are added to the screw as so-called dynamic mixing elements.
Despite a long history of the field, the great variety of these mixing elements reported, e.g., by Gale \cite{gale2009mixing}, or Campbell and Spalding \cite{Campbell2013analyzing} indicate that designing an optimal mixing element shape is still an unsolved task.
A good overview of recent works towards assessing the various designs and current approaches to create improved mixing elements can be found in \cite{hopmann2020method}.\par
In recent years, the design of mixing elements is increasingly aided by numerical simulation. 
Potente and T\"obben, e.g., investigate different designs of the Maddock mixer using unwound screw channels \cite{Potente2002Maddock}.
Over 5000 finite-element simulations are conducted to derive a functional that relates the pressure-throughput behavior of the studied mixers to specific design quantities like the helix angle.
The study can be considered a first step in the direction of reduced-order modeling.
Another example is the work by Celik et al. who -- based on simulations -- assess distributive mixing based on particle-tracking and dispersive mixing based on particle deformations \cite{Celik2017Einfaerbung}.
%However, their analysis focuses on dispersive mixing only, and again, different designs are compared, but no optimization occurs.
Finite-element simulations are also used by Marschik et al. \cite{Marschik2018}, as well as Roland et al. \cite{Roland2019FEM} to investigate the isothermal flow field around block-head mixing screws and pineapple mixers, respectively.
Both works also conduct computational design studies of the analyzed mixers.\par
In the literature listed above, the focus lies on the comparison of different, yet predefined mixing-element designs. 
Utilizing numerical methods to perform an actual design optimization is still out of reach. 
Instead, design optimization is so far primarily performed based on physical experiments.\par
Conversely, in many other areas, numerical design based on finite-element simulations is already very mature.
Such simulation-based design frameworks -- given that suitable initial designs are often available -- typically rely on shape optimization.
One essential aspect of shape optimization is finding suitable parameterizations of  the computational domain.
Among the common choices are radial basis functions \cite{RBF_airfoil, Kobbelt2004}, neural networks \cite{VAE_shapeMorph}, and free-form deformation \cite{RozzaFFD2010, Rozza_FFD13}.
Also very widely used are parameterizations based on the spline-definition within CAD models.
Using a linear-elasticity-based mesh update method to connect the boundary shape modifications to updates of the volumetric computational mesh, this technique is used by Elgeti et al. in the context of shape optimization of extrusion dies \cite{Elgeti12}.
While these approaches are employed to find a low-dimensional parameterization, Hojjat et al., in contrast, focus on an extremely versatile shape parameterization using vertex-based shape modifications \cite{Bletzinger_vertex_method}.
Due to the huge number of optimization variables, the vertex-based approach, however, relies upon the high efficiency of adjoint methods for shape optimization.\par
In view of these advances in numerical design, this study's purpose is bridging the still existing gap between analysis and design in the field of mixing elements applying current methods of numerical shape optimization.
In particular, we aim to shape-optimize an individual mixing element rather than a complete mixing section.
The design is based on a representative subproblem: Instead of considering the flow through the entire mixing section, the optimization is based on a single mixing element.
The paper is structured as follows:
First, we review the general problem statement and introduce a novel objective formulation for distributive mixing in Sec.~\ref{sec:shpopt_framework}.
In Sec.~\ref{sec:shape_param}, we review Free-Form Deformation as a universal shape parameterization tool, before we outline the governing equations of the investigated flow and mesh-update problems in Sec.~\ref{sec:goveq}.
The specific geometrical set-up is introduced in Sec.~\ref{sec:model}, leading up to the application case and corresponding experimental validation in Sec.~\ref{sec:results}.
Finally, we conclude with a discussion of further steps in Sec.~\ref{sec:conclusion}.

\section{Shape optimization problem}
\label{sec:shpopt_framework}
To derive a computational model, we first review the general problem of designing mixing sections, which is: \\
\textit{Find the optimal form, distribution, and number of mixing elements such that mixing is enhanced.} \\
This problem can be mathematically formulated as a Partial Differential Equation (PDE)-constrained linear programming problem, which takes the following form:
\begin{subequations}\label{eq:template_opt_prob}
  \begin{alignat}{3}
    &\min_{\boldsymbol{\alpha}}{J}\left(\boldsymbol{\alpha}\right)\\
    \text{s.t.}\hspace{0.5em}		     &\mathbf{F}\left(\boldsymbol{\alpha}\right) = \mathbf{0} \qquad  &&\text{in} \hspace{0.5em} \Omega\left(\boldsymbol{\alpha}\right), \\
					     &\alpha_i \geq \alpha_{min,i}, &&i=1,...n_{\alpha}, \label{eq:l_bounds}\\
					     &\alpha_i \leq \alpha_{max,i}, &&i=1,...n_{\alpha}. \label{eq:u_bounds}
  \end{alignat}
\end{subequations}
Here, $J$ denotes the objective function -- in this scenario describing the mixing quality --, $\boldsymbol{\alpha}$ the optimization parameters -- here describing the form and distribution of mixing elements --, and $\mathbf{F}$ the PDE constraints defined on the computational domain $\Omega$ -- here, the PDE describes the flow solution inside the extruder.
Furthermore, Eqs.~\eqref{eq:l_bounds}  and \eqref{eq:u_bounds}  allow for bound constraints on the optimization parameters.
%where the latter two conditions describe bound constraints on the optimization variables $\boldsymbol{\alpha}$, and %$\mathbf{F}\left(\boldsymbol{\alpha}\right)=\mathbf{0}$ denotes the set of governing partial differential equations that needs to be fulfilled on the
%computational domain $\Omega\left(\boldsymbol{\alpha}\right)$.
Solving \eqref{eq:template_opt_prob} yields the optimal mixing section's shape with respect to the chosen objective function.
Due to the complexity of the considered problem, we restrict the optimization to an appropriate submodel, which will be described in the next section along with an illustration of the choices for $J$, $\boldsymbol{\alpha}$, and $\mathbf{F}$.

\subsection{Submodel choice}
\label{subsec_submodel}
In an ideal situation, one would try to solve the general problem as stated above, thereby optimizing the complete mixing section consisting of several mixing elements.
The problem with this approach is the curse of dimensionality.
To simultaneously enable a varying number and position of all mixing elements, a very detailed shape parameterization is required.\par
These drawbacks motivate a split of the shape optimization problem, which is outlined by Eusterholz and Elgeti in \cite{eusterholz2018esaform} and reviewed in the following:
Instead of solving one global shape optimization problem, we propose to solve two smaller problems sequentially:
(1) Shape optimize a single mixing element using standard spline-based methods to determine an optimized shape, and (2) distribute these mixing elements on the screw body subsequently.
With this, the final setup is one linear programming and one mixed integer programming problem.
While both problems, being PDE-constrained, remain computationally expensive, the proposed split allows separating the many variables needed to shape-optimize a mixing element from the integer variables determining the number of elements.
This split reduces the computational complexity of the optimization problem and allows to set up two simulations, each tailored to the different optimizations, yielding further performance gains.
In this paper, we focus on the first, i.e., the mixing-element optimization, and demonstrate how the proposed framework already increases the extruder's performance.

\subsection{Objective function for distributive mixing}
\label{subsec:obj}
Whether or not a screw performs well will in this work be measured by it's ability to enhance distributive mixing.
To quantify mixing, scalar measures $J$ need to be found.
While the actual implementation will be described in Sec.~\ref{sec:goveq}, the proposed measure's physical motivation is described in the following.\par
%We thus aim to homogenize the concentration of one material mixed into another one.
%Typical measures for distributive mixing usually quantify to what extend two -- possibly immiscible -- phases become indistinguishable.
Typical measures for distributive mixing used in simulation often rely on tracking massless particles convected with the flow.
Using the traces of these particles, distributive mixing is often evaluated based on residence time distribution or the concentration's standard deviation, as summarized in \cite{Campbell2013analyzing}.
Other measures inspired by Danckwerts'early work \cite{Danckwerts1952} utilize pairwise correlation functions \cite{YangCorrelation09, Wang01, WangZloczower03}.
An improvement to such correlation functions was proposed by the group of Manas-Zloczower introducing Renyi entropies \cite{Wang01, WangZloczower03}.
While all these measures have proven to be suitable to assess mixing in well-mixed systems, we faced difficulties  using such mixing indices to measure only subtle mixing increases.
Since in our submodel, mixing is expected to be only slightly enhanced -- as a consequence of investigating the influence of only a single mixing element -- we revisit experimental evidence to propose a different measure.
One common mixing measure in experiments is the increase in interfacial area.
Using the interfacial area as a mixing measure is motivated by the observation that alternating layers of two materials become indistinguishable as soon as their so-called striation thickness approaches zero.
As found in various studies, \cite{Mohr57, Spencer51, gale2009mixing, klie2015phd, covas1995rheological} the striation thickness, or its inverse the interfacial area, is a well-known and proven quantifier for distributive laminar mixing \cite{Mohr57}.
Starting from the experimental reliability of this parameter, we derive an estimator that allows us to quantify the interface enlargement for moderate flow reorientations.\par
To estimate the interface enlargement numerically, we define rectangular areas at the inflow of the computational domain, track them towards the outflow, and compare their perimeter length before and after passing a mixing element.
It should be noted that from a fluid mechanics point of view, in this approach we determine the  increase of the cross-section perimeter length between the inlet and the outlet from the computational domain of different stream tubes. 
A distortion of the stream tube leads to larger perimeter increases and enlarged interfacial areas after flowing around a mixing element.
A detailed description of the implementation and numerical approximations of this concept is given in Sec.~\ref{subsec:goveq_mixing}.

\section{Shape parameterization}
\label{sec:shape_param}
\begin{figure}[!htb]
	\centering
	\input{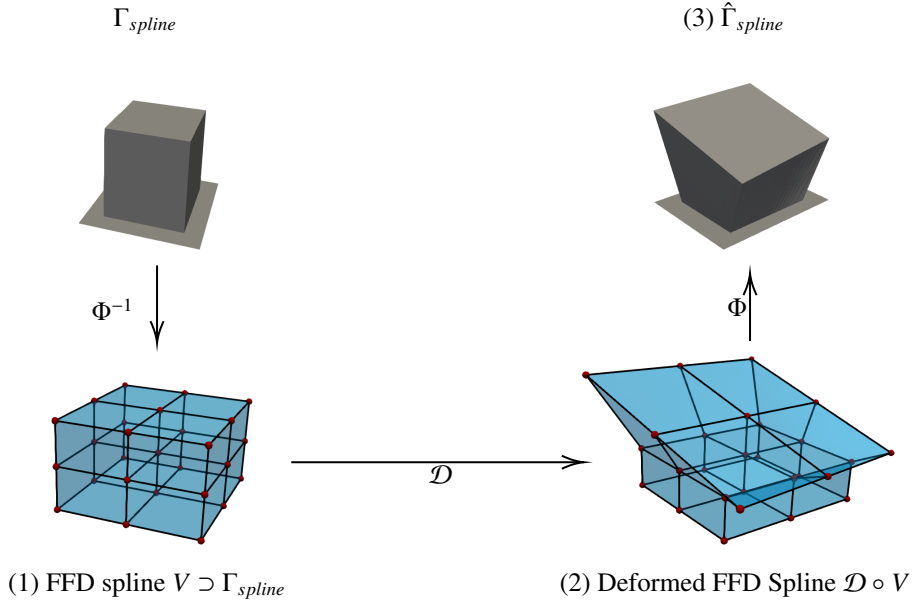}
	\caption{
		Complete FFD workflow in three steps:
		(1) construction of an enclosing domain $V\supset\Gamma_{spline}$ of the mixing element $\Gamma_{spline}$ used as the control volume of the FFD spline,
		(2) deformation of the FFD spline, and
		(3) construction of the deformed shape through evaluation of the deformed FFD spline for all points $\mathbf{x} \in \Gamma_{spline}$ at their corresponding parametric coordinates in $V$, $\Phi^{-1}\left(\mathbf{x}\right)$
		.
	}
	\label{fig:ffd1}
\end{figure}

The shape parameterization, i.e., the choice of optimization variables controlling the shape deformation, is essential for an effective yet fast shape optimization.
One possible shape parameterization is given by the finite element mesh itself.
However, the resulting large number of optimization variables requires sophisticated adjoint solvers \cite{Bletzinger_vertex_method}.
Due to given constraints on the utilized code, but primarily motivated by the use of commercial simulation software in the industry, our approach relies on "black-box" optimization techniques, where we do not require an analytical evaluation of the Jacobian of our optimization problem \eqref{eq:template_opt_prob}.
%, especially not of the term $\nabla_{\boldsymbol{\alpha}}\mathbf{F}\left(\boldsymbol{\alpha}\right)$ -- the gradient of the PDE-constraints.
While  many techniques exist to solve such black-box optimization problems, a standard limitation is that the number of optimization variables has to be kept small.
One way to ensure such a small number of optimization variables is using CAD shape representations.
In using such CAD data, the key idea is to utilize the already existing geometry description based on non-uniform rational B-Splines (NURBS).
These splines are based on a polynomial interpolation of typically only a few control points, which may be directly used as optimization variables \cite{Elgeti12}.
For a thorough description of NURBS's theory, the reader is referred to the book of Piegl and Tiller \cite{piegl1996nurbs}.
With regard to the optimization of complete mixing sections it should be noted that although using CAD splines works well for shape optimizations where only minor shape changes are expected, it is in general not suitable for the shape optimization of mixing sections.
The reason why mixing sections need a more versatile shape parameterization is that the number and position of the mixing parts should be variable.
While CAD splines generally allow the deformation of existing geometry features, they do not offer the required variability to add new features, i.e., new mixing elements.
Adding and subtracting mixing elements, however, may be achieved using trimmed NURBS, i.e., Boolean additions of multiple splines.
Unfortunately, Boolean operations leads to more complex Mixed Integer Programming problems.
In this work, however, we only focus on the single-mixing-element shape optimization as described in Sec.~\ref{subsec_submodel}.
Thus, a first reduction of the number of shape parameters can be achieved by limiting the parameterization to the surface representation of the mixing element; this can be viewed as an extension of Elgeti et al. \cite{Elgeti12}.
However, even with this restriction to surface representations, the underlying NURBS models may still have to be reworked in order to ensure an appropriate, yet efficient enough resolution of the surface parameterization.
In order to avoid this step, we resort to free-form deformation (FFD) \cite{Sederberg1986}; in this approach, the surface is embedded into a volume spline.
Due to the individualizable choice of the number of control points of the volume spline, the number of optimization parameters can be selected by the user while at the same time, all deformations are guaranteed to be smooth.  
FFD is a shape deformation technique that is used in various fields to allow for either (a) a large parameter reduction or (b) smooth deformations.\\
The general three-step FFD workflow is depicted in Fig.~\ref{fig:ffd1}.
First, the object under consideration is embedded into a volume spline of arbitrary dimension; meaning that every point coordinate $\mathbf{x}$ on the object is mapped to a specific parametric coordinate in the volume spline, denoted as $\Phi^{-1}\left(\mathbf{x}\right)$. 
During the shape optimization process, the control points of the volume spline can be displaced, leading to a continuous displacement field which finally smoothly maps every initial surface point to a new, deformed configuration.
This displacement is denoted via the deformation operator $\mathcal{D}$. 
%A deformed geometry is obtained by moving only a few control points, given any interpolation-based shape representation, and finally imposing this mapping onto the support points (3).
The resulting coordinates of the deformed object can be computed as
\begin{equation}
	\hat{\mathbf{x}}=\mathcal{D}\circ\Phi\left(\Phi^{-1}\left(\mathbf{x}\right)\right).
	\label{eq:}
\end{equation}
Note that in the general case, $\Phi^{-1} \hspace{0.15em}:\hspace{0.15em} \mathbf{x} \mapsto \boldsymbol{\xi}$ can be obtained using Newton's method solving $\mathbf{x} = \Phi\left(\boldsymbol{\xi}\right)$.
%However, in FFD, it is more common to employ volume splines, where the parameter space and the physical space have the same dimension. Then, a more convenient linear mapping $[\mathbf{x}_{min}, \mathbf{x}_{max}] \rightarrow [\boldsymbol{\xi}_{min}, \boldsymbol{\xi}_{max}]$ can be employed given suitable splines \cite{Rozza_FFD13}.
Due to its aforementioned properties, FFD has recently gained much attention in shape optimization, often combined with model order reduction techniques \cite{RozzaFFD2010, Rozza_FFD13,  salmoiraghi2018free}.
In our case, we primarily utilize FFD to provide a robust interface to arbitrary CAD splines.
The final and necessary reduction in parameters is achieved by defining a set of admissible degrees of freedom in advance that relates the movement of control point groups.\par
In the presented example, the complete shape parameterization is obtained by constructing a Cartesian spline around the CAD geometry of one mixing element and allowing the following six deformations applied to the embedding spline:
(1--3) rotation  of the spline's top plane around $x$, $y$, and $z$, (4) change in height, (5) radial expansion of the top to allow for conical shapes, and finally (6) rotation of the complete spline around its $z$ axis.

\section{Governing equations}
\label{sec:goveq}
In this section, we describe the governing equations of the computational model. 
Next to the flow equations, this includes a PDE for mesh adaptation, and a transport equation to mimic particle tracking.
Furthermore, we outline the discretization approach. 
\subsection{Flow equations}
The governing equations for the computational model are the steady incompressible Navier-Stokes equations: 
\begin{subequations}
	\begin{alignat}{3}
		\rho\left(\mathbf{u}\cdot \nabla\right)\mathbf{u} - \nabla \cdot \boldsymbol{\sigma}&=\mathbf{0} \qquad && \text{in} \enskip \Omega, \\
		\nabla \cdot \mathbf{u} &=0 \qquad &&  \text{in} \enskip \Omega ,\\
		\mathbf{u} &= \mathbf{u}_D \qquad && \text{on} \enskip \Gamma_{D_{\mathbf{u}}} ,\\
		\mathbf{n}\cdot\boldsymbol{\sigma} &= \mathbf{h}_{\mathbf{u}} \qquad && \text{on} \enskip \Gamma_{N_{\mathbf{u}}},
	\end{alignat}
\end{subequations}
combined with a steady heat equation
\begin{subequations}
	\begin{alignat}{3}
		\rho c_p\mathbf{u}\cdot \nabla T -  \lambda \Delta T +\boldsymbol{\sigma}\colon \nabla \mathbf{u} &= 0 \qquad && \text{in} \enskip \Omega ,\\
		T &= T_D \qquad &&  \text{on} \enskip \Gamma_{D_T}, \\
		\mathbf{q} &= \mathbf{h}_T \qquad && \text{on} \enskip \Gamma_{N_T}.
	\end{alignat}
\end{subequations}
Here, $\mathbf{q}=-\lambda \nabla T$ denotes the heat flux, and $\Gamma_{D_\Box}$ and $\Gamma_{N_\Box}$ denote the Dirichlet and Neumann part of the boundary of the domain  $\Omega$ for the respective variables $\mathbf{u}$ -- the flow field -- and $T$ -- the temperature -- represented by the placeholder $\Box{}$.
Using the same placeholder, $\mathbf{h}_{\Box{}}$ denotes  the respective Neumann boundary fluxes.
The material parameters density, specific heat, and thermal conductivity are given as $\rho$, $c_p$, and $\lambda$.
Finally using the rate-of-strain tensor, ${\bm{\varepsilon}}\left(\Box{}\right) = \frac{1}{2}\left(\nabla \Box{}+\nabla \Box{}^T\right)$, and the pressure $p$, the stress tensor for Newtonian fluids, $\boldsymbol{\sigma}$, is modeled as:
\begin{equation}
	\boldsymbol{\sigma} = -p\mathbf{I} + 2 \eta  \varepsilon\left(\mathbf{u}\right).
\end{equation}
Both equation systems are coupled through viscous dissipation and the constitutive model for the viscosity $\eta$ which depends on the shear rate $\dot{\gamma} = \sqrt{2 \varepsilon\left(\mathbf{u}\right)\colon\varepsilon\left(\mathbf{u}\right)}$.
%Even though the flow regime admits describing the flow by the linear Stokes equations, the material model introduces another non-linearity to the system of fluid equations, and thus we find it beneficial to solve the non-simplified equations.
As a constitutive model for the viscosity, we consider the three-parameter  Carreau model \cite{carreau72} which, using the notation of  \cite{michaeli2009}, can be written as:
\begin{equation}
	\eta\left(\dot{\gamma}\right) = \frac{A}{\left(1+B\dot{\gamma}\right)^C}\enskip,
	\label{eq:carreau}
\end{equation}
where $A$ denotes the viscosity at zero shear rate, $\frac{1}{B}$ indicates the transition to the shear-thinning regime and $C$ describes the rate at which the viscosity decreases.
%For the given material we prescribe $A=\SI{9,472}{\pascal \second}$, $B=\SI{0.1871}{\second}$, and $C=0.655$.
While pressure corrections for the viscosity are not taken into account, a Williams-Landel-Ferry temperature correction is employed.
Following the notation of Osswald and Rudolph \cite{PolymerRheologyCh3}, and noting that $A=\eta_0$, a correction for $A$ is given as:
\begin{equation}
	\log\frac{A\left(T\right)}{A\left(T^*\right)}=\frac{8.86\left(T^*-T_s\right)}{101.6+T^*-T_s}-\frac{8.86\left(T-T_s\right)}{101.6+T-T_s}\enskip.
	\label{eq:}
\end{equation}
In this temperature correction, $T^*$ denotes the measurement temperature at which material parameters are determined and $T_s$ denotes a reference temperature often chosen as a constant or pressure-dependent offset to the glass transition temperature.

\subsection{Mesh adaptation}
In Sec.~\ref{sec:shape_param}, we introduced our approach to FFD, which is applied to the CAD surface representation.
This entails a total of three steps until a volume mesh suitable for finite-element flow computations is obtained: (1) update the CAD representation using FFD, (2) adapt the spline surface mesh, $\Gamma_{spline}$, accordingly by reevaluating the spline at the parametric coordinates connected to the FEM nodes, and (3) adjust the volume mesh to this deformation in order to enhance mesh quality.
Step (3) is performed using the Elastic Mesh Update Method (EMUM) as described by Johnson and Tezduyar \cite{Johnson94a}. In EMUM, the volume mesh is treated as an elastic body reacting to displacements imposed at its boundaries.
These displacements are exactly the displacements of the surface mesh (cf. Step (2)).
EMUM's governing equations, with the displacement $\mathbf{z}$ being the unknown, are given in Eqs.~\eqref{eq:EmumGovEq1} to \eqref{eq:EmumGovEq3}:\\
\begin{align}
	\label{eq:EmumGovEq1}
	\nabla \cdot \mbox{${\bm{\sigma}}$}_{\mathrm{mesh}} &= {\mathbf 0} \,, \\
	\mbox{${\bm{\sigma}}$}_{\mathrm{mesh}}({\mathbf z}) &= \lambda_{\mathrm{mesh}} \left(
	{\mathrm{tr}} \, {\bm{\varepsilon}} ({\mathbf z})\right){\mathbf
		I} + 2 \mu_{\mathrm{mesh}} {\bm{\varepsilon}}({\mathbf z}) \,.% \\
%	{\bm{\varepsilon}}_{\mathrm{mesh}} ({\mathbf z}) &= 
%	\text{sym}\left(\nabla \mathbf{z}\right) \,.
	%\frac{1}{2} \left(\nabla
	%{\mathbf z} + (\nabla {\mathbf z})^T\right) \,.
	\label{eq:EmumGovEq3}
\end{align}
Note that $\lambda_{\mathrm{mesh}}$ and $\mu_{\mathrm{mesh}}$ are no longer used in their original sense as Lam\'e parameters, but they can be employed as free parameters that influence the mesh quality. 
To achieve optimal mesh relaxation, the boundary is partitioned into three parts as $\Gamma = \Gamma_{spline} \cup \Gamma_{fixed} \cup \Gamma_{slip}$ as shown in Fig.~\ref{fig:simDom}.
Tangential movement is allowed on $\Gamma_{slip}$, whereas we restrict any movement on $\Gamma_{fixed}$, and prescribe the displacement obtained from the spline deformation on $\Gamma_{spline}$.
Using this approach, we are able to use boundary-conforming meshes without the need for remeshing.

\subsection{Mixing assesment using advection-based particle tracking}
\label{subsec:goveq_mixing}
In Sec.~\ref{subsec:obj}, the novel objective formulation based on the deformation of stream tubes during the mixing is introduced.
Based on this knowledge, in this section, we elaborate on how such an objective formulation can be implemented.\par
In general, stream tubes can be approximated via a discrete set of streamlines, which in turn can be determined by tracking a set of massless particles.
However, the integration of particle paths is prone to errors, especially in the presence of no-slip boundary conditions.
Therefore, instead of tracking individual particles, we consider each particle as a scalar value that is advected through the flow field.
In the following, we will explain how we set up the advection problems to study the movement of single particles.\par
First, we define a sampling region \(\Gamma_s\) as a portion of the inflow boundary \(\Gamma_{in}\) located upstream of the mixing element.
This choice is depicted in Fig.~\ref{fig:simDom} as the dotted area and marks the starting point of the evaluated stream lines.
We restrict the sampling area to this area of the inflow to limit the computational demand.
We choose this specific placement since we expect most mixing to happen in the vicinity of the mixing element.
Within the sampling area, we refine the mesh such that the support of each shape function is reduced.
To mimic a particle tracking approach, we then solve a  set of steady advection equations for the scalar property \(c_j\):
\begin{alignat}{3}
	\mathbf{u}\cdot \nabla c_j &= 0 &&\text{in} &&\Omega, \\
	c_j                            &= 0  \quad     &&\text{on}       \quad \quad             &&\Gamma_{in} \setminus \Gamma_s, \\
	c_j\left(\mathbf{x}_i\right) &= \delta_{ij}  &&\forall \mathbf{x}_i \in     &&\Gamma_{s},
\end{alignat}
where $j$ indicates an iteration over each mesh node located inside the sampling area $\Gamma_s$.
This approach leads to a set of $n$ distributions of \(c_j, j\in[1,n]\) at the outflow domain, with $n$ being the number of nodes in $\Gamma_s$.
We make the assumption that given a sufficient resolution of \(\Gamma_{out}\), regardless of numerical diffusion, the maxima of \(c_j\) at \(\Gamma_{s}\) and \(\Gamma_{out}\) correspond to the same streamline.
This notion of streamlines allows correlating particle positions between in- and outflow as a set of particle movement vectors.
Projecting these vectors onto the outflow boundary plane, $\Gamma_{out}$, one arrives at a set of two-dimensional displacement vectors $\mathbf{d}_j\colon \mathbf{x}\left(\max{c_j}\right)\big|_{\Gamma_{in}}\mapsto \mathbf{x}\left(\max{c_j}\right)\big|_{\Gamma_{out}}$ that defines an inflow to outflow mapping for each particle.\par
In the limit of vanishing element sizes, this approach constructs a linear mapping $A$, such that
\begin{equation}
	\mathbf{c}\left(\mathbf{x}_{out}\right) = A \mathbf{c}\left(\mathbf{x}_{in}\right).
	\label{eq:linMap}
\end{equation}
In theory, using Eq.~\eqref{eq:linMap} enables applying correlation techniques known from image processing.
Such methods, in general, quantify the correlation between inflow and outflow, but in addition allow investigating at which length scales mixing is enhanced.
In the discretized framework, however, it is found that the condition that any two streamlines must not cross each other can be violated, especially when only evaluating the positions of maximal concentration  $\mathbf{x}\left(\max{c_j}\right)$.\par
To instead quantify laminar mixing based on the displacement vectors \(\mathbf{d}\), several options exist.
Primitive choices like averaged $L_1$ or $L_2$ norm cannot predict mixing properly, especially since such measures indicate increased mixing even in cases where the melt is subjected to a uniform translational displacement.
While a sensitivity against this can be incorporated by evaluating $\rm{var}\left(\mathbf{d}\right)$, we found that no combination of any of these measures is globally monotonous and sensitive at the same time.
Therefore, we choose to quantify mixing by approximating the increase of the interfacial area between two virtual phases.\par
As discussed before, split \(\Gamma_s\) into \(m \times n\) rectangular non-overlapping subdomains which resemble stream tubes.
The goal, again, is to compare the stream tube's perimeter between inflow and outflow as a measure of the interfacial area  increase. 
The approach for a total of $k$ subdomains is as follows:
Let \(X_k\) be the set of mesh nodes \(\mathbf{x}_k\) in the $k^{\text{th}}$ subdomain, \(D_k\) the set of displacement vectors \(\mathbf{d}\) for all nodes in \(X_k\), \(\mathcal{C}\) an operator to construct the convex hull around a set of points, and \(\mathcal{L}\) a measure
for the length of a curve.
Then we can approximate the increase in interfacial area as the sum of the increase of all convex hull lengths as:
\begin{equation}
	J = \sum_{k=1}^{nm} \mathcal{L}\left(\mathcal{C}\left(X_k\right)\right) - \mathcal{L}\left(\mathcal{C}\left(X_k+D_k\right)\right).
	\label{eq:obj}
\end{equation}
An illustration is given in Fig.~\ref{fig:three_objectives} showing the convex hulls at the outflow of $2 \times 2$ initially rectangular subdomains.
The restriction to convex hulls here is motivated by the fact that a unique convex hull exists for each set of points, and so do computational frameworks to identify them \cite{SciPy2020, Barber96thequickhull}.
It should be noted that the proposed method, though exemplified on a box domain, can not only be applied to arbitrary domains with coplanar in- and outflows but is without modification also applicable to experiments.

\subsection{Flow solver}
\begin{figure}[!htb]
	\centering
	\tikzset{every picture/.style={line width=0.75pt}} %set default line width to 0.75pt

\begin{tikzpicture}[x=0.75pt,y=0.75pt,yscale=-1,xscale=1]
%uncomment if require: \path (0,310); %set diagram left start at 0, and has height of 310

%Straight Lines [id:da08701824394621938]
\draw    (130.5,68) -- (130.5,86) ;
\draw [shift={(130.5,88)}, rotate = 270] [color={rgb, 255:red, 0; green, 0; blue, 0 }  ][line width=0.75]    (10.93,-3.29) .. controls (6.95,-1.4) and (3.31,-0.3) .. (0,0) .. controls (3.31,0.3) and (6.95,1.4) .. (10.93,3.29)   ;
%Shape: Rectangle [id:dp2482477698576464]
\draw   (90.5,35) -- (504.5,35) -- (504.5,211) -- (90.5,211) -- cycle ;
%Straight Lines [id:da5554072046523156]
\draw    (130.5,113) -- (130.5,131) ;
\draw [shift={(130.5,133)}, rotate = 270] [color={rgb, 255:red, 0; green, 0; blue, 0 }  ][line width=0.75]    (10.93,-3.29) .. controls (6.95,-1.4) and (3.31,-0.3) .. (0,0) .. controls (3.31,0.3) and (6.95,1.4) .. (10.93,3.29)   ;
%Straight Lines [id:da6209936830955936]
\draw    (130.5,158) -- (130.5,176) ;
\draw [shift={(130.5,178)}, rotate = 270] [color={rgb, 255:red, 0; green, 0; blue, 0 }  ][line width=0.75]    (10.93,-3.29) .. controls (6.95,-1.4) and (3.31,-0.3) .. (0,0) .. controls (3.31,0.3) and (6.95,1.4) .. (10.93,3.29)   ;

% Text Node
\draw    (113,41) -- (152,41) -- (152,67) -- (113,67) -- cycle  ;
\draw (120,49) node [anchor=north west][inner sep=0.75pt]   [align=center] {FFD};
% Text Node
\draw    (104,87) -- (161,87) -- (161,113) -- (104,113) -- cycle  ;
\draw (111,95) node [anchor=north west][inner sep=0.75pt]   [align=center] {EMUM};
% Text Node
\draw    (104,132) -- (161,132) -- (161,158) -- (104,158) -- cycle  ;
\draw (113,140) node [anchor=north west][inner sep=0.75pt]   [align=center] {FLOW};
% Text Node
\draw    (112,177) -- (153,177) -- (153,203) -- (112,203) -- cycle  ;
\draw (119,184) node [anchor=north west][inner sep=0.75pt]   [align=center] {OPT};
% Text Node
\draw (175,46) node [anchor=north west][inner sep=0.75pt]   [align=left] {spline based shape parameterization};
% Text Node
\draw (175,91) node [anchor=north west][inner sep=0.75pt]   [align=left] {domain adaptation based on linear elasticity};
% Text Node
\draw (176,137) node [anchor=north west][inner sep=0.75pt]   [align=left] {forward simulation for flow and particle tracking};
% Text Node
\draw (176,182) node [anchor=north west][inner sep=0.75pt]   [align=left] {optimization algorithm};

\end{tikzpicture}
	\caption{Optimization framework including mesh update modules, simulation and optimization driver.}
	\label{fig:framework}
\end{figure}
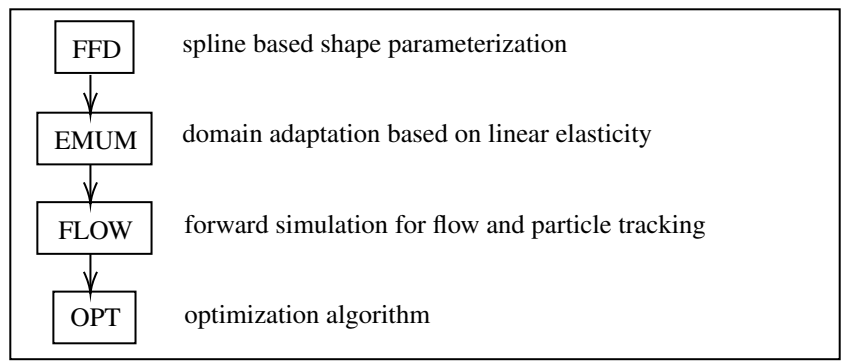

 All the presented governing equations for flow,  temperature, and concentration fields are discretized using Newton's method for linearization and Galerkin-Least-Squares (GLS)-stabilized linear finite elements.
The resulting equation systems are solved iteratively using a preconditioned GMRES algorithm, implemented into the highly-parallel in-house flow solver. To account for the shear-thinning behavior, we internally couple the flow and heat equations in Gauss-Seidel fashion.
Since these concentrations $c_j$ (cf. Sec.~\ref{subsec:goveq_mixing}) do not affect the flow field, we solve those sequentially in a weakly coupled fashion after converging the flow-temperature system.
The complete framework consists of the building blocks shown in Fig.~\ref{fig:framework}.
\section{Test case generation}
{\label{sec:model}
In this section, we will introduce the chosen test case.
In particular, we explain how this test case is derived as a sub part of the unwound screw channel.
Finally, we detail on the utilized boundary conditions and give a list of the chosen material parameters.\par
In most screws, the flow channel features a regular pattern of several mixing elements.
We make use of this pattern and segment the flow channel into non-overlapping rectangular domains such that each domain covers exactly one mixing element.
These rectangular domains are the basis for our test case.
Since due to the pattern, all such domains are equal, we only consider one.\par
As explained, the domain is constructed around one mixing element whose orientation is discussed in the following.
In contrast to the box domain, the mixing element is not aligned with the $y$-axis, but rotated in negative $z$ direction by \(5.7\degree\).
This $z$-rotation models the mixing element's angle of attack.
This angle approximates the streamlines' pitch around the screws rotational axis in an actual extruder.
The pitch angle's small magnitude -- relating circumferential movement to actual throughput -- reveals that the flow inside the extruder mainly follows the screw's rotational movement.
Consequently, we also neglect non-circumferential flow in our model.\par
The dimensions and general arrangement of our computational model are depicted in Fig~\ref{fig:simDom}, where \(x\) corresponds to the extruder's rotation axis.
As the computational domain resembles a part of the unwound screw, the height matches the flow channel height of an actual \(60\unit{mm}\) screw, and the width resembles the distance between adjacent mixing elements.
The domain length $L$ is chosen empirically such that no backward flow over $\Gamma_{out}$ into the domain occurs and satisfies $\frac{L}{D} \approx 0.4 $.

%As motivated in Sec.~\ref{subsec_submodel}, we simulate the flow around single mixing element in a part of the unwound screw channel.
%Fig.~\ref{fig:simDom} depicts this part of the unwound screw and illustrates the computational domain's dimension.
%While the purpose is only to provide insight into the flow around a single mixing element, these dimensions are inspired by an actual screw, i.e., the height matches the flow channel height of an actual \(60\unit{mm}\) screw and the width resembles the
%distance between adjacent mixing elements in an unwound screw.
%The domain length $L$, in contrast, is chosen empirically such that no backward flow over $\Gamma_{out}$ into the domain occurrs and measures $\frac{L}{D} \approx 0.4 $.
%It should be noted that the orientation is chosen such that \(x\) corresponds to the extruder's rotation axis.
%This orientation is the result of fluid simulations of the complete screw including the rotational movement which show that the flow is mainly following the screw in circumferential direction.
%Experimental investigations imply further that in fact the flow's pitch angle is only \(5.7\degree\).
%This pitch angle is incorporated in the simulation by a rotation of the mixing element in negative $z$-direction.
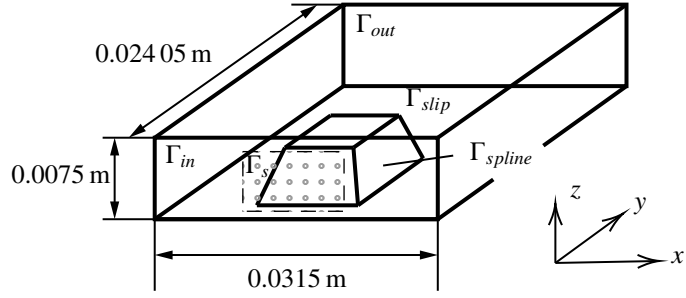
\begin{figure}[!htb]
  \centering
  \tikzset{
pattern size/.store in=\mcSize,
pattern size = 5pt,
pattern thickness/.store in=\mcThickness,
pattern thickness = 1.0pt,
pattern radius/.store in=\mcRadius,
pattern radius = 1pt}
\makeatletter
\pgfutil@ifundefined{pgf@pattern@name@_ax055x350}{
\makeatletter
\pgfdeclarepatternformonly[\mcRadius,\mcThickness,\mcSize]{_ax055x350}
{\pgfpoint{-0.5*\mcSize}{-0.5*\mcSize}}
{\pgfpoint{0.5*\mcSize}{0.5*\mcSize}}
{\pgfpoint{\mcSize}{\mcSize}}
{
\pgfsetcolor{\tikz@pattern@color}
\pgfsetlinewidth{\mcThickness}
\pgfpathcircle\pgfpointorigin{\mcRadius}
\pgfusepath{stroke}
}}
\makeatother
\tikzset{every picture/.style={line width=0.75pt}} %set default line width to 0.75pt
\begin{tikzpicture}[x=0.75pt,y=0.75pt,yscale=-1,xscale=1]
%uncomment if require: \path (0,300); %set diagram left start at 0, and has height of 300

%Shape: Rectangle [id:dp8695735038286327]
\draw  [line width=1.5]  (179,194) -- (320.3,194) -- (320.3,235) -- (179,235) -- cycle ;
%Straight Lines [id:da28378117790533053]
\draw [line width=1.5]    (273.2,127.05) -- (179,194) ;
%Straight Lines [id:da5615496570179084]
\draw [line width=1.5]    (273.2,168.05) -- (179,235) ;
%Straight Lines [id:da518845160801926]
\draw [line width=1.5]    (414.5,127.05) -- (320.3,194) ;
%Straight Lines [id:da5705336604419526]
\draw [line width=0.75]    (257.75,128.79) -- (168.45,192.26) ;
\draw [shift={(166,194)}, rotate = 324.6] [fill={rgb, 255:red, 0; green, 0; blue, 0 }  ][line width=0.08]  [draw opacity=0] (8.93,-2.29) -- (0,0) -- (8.93,2.29) -- cycle    ;
\draw [shift={(260.2,127.05)}, rotate = 144.6] [fill={rgb, 255:red, 0; green, 0; blue, 0 }  ][line width=0.08]  [draw opacity=0] (8.93,-2.29) -- (0,0) -- (8.93,2.29) -- cycle    ;
%Shape: Rectangle [id:dp05088925318353121]
\draw  [line width=1.5]  (273.2,127.05) -- (414.5,127.05) -- (414.5,168.05) -- (273.2,168.05) -- cycle ;
%Straight Lines [id:da6560330089410432]
\draw [line width=1.5]    (320.3,235) -- (179,235) ;
%Straight Lines [id:da5374059432809515]
\draw [fill={rgb, 255:red, 0; green, 0; blue, 0 }  ,fill opacity=1 ][line width=1.5]    (266.8,182.8) -- (244,199) ;
%Straight Lines [id:da6825269545935158]
\draw [fill={rgb, 255:red, 0; green, 0; blue, 0 }  ,fill opacity=1 ][line width=1.5]    (301,182.8) -- (278.2,199) ;
%Straight Lines [id:da3633541419106372]
\draw [fill={rgb, 255:red, 0; green, 0; blue, 0 }  ,fill opacity=1 ][line width=1.5]    (278.2,199) -- (244,199) ;
%Straight Lines [id:da860924902224956]
\draw [fill={rgb, 255:red, 0; green, 0; blue, 0 }  ,fill opacity=1 ][line width=1.5]    (301,182.8) -- (266.8,182.8) ;

%Straight Lines [id:da3922046037567153]
\draw [line width=1.5]    (244,199) -- (230,228) ;
%Straight Lines [id:da9535258760167132]
\draw [line width=1.5]    (278.2,199) -- (280.49,228) ;
%Straight Lines [id:da2502148496462797]
\draw [line width=1.5]    (301,182.8) -- (313,204.4) ;
%Straight Lines [id:da6775260688023841]
\draw [fill={rgb, 255:red, 0; green, 0; blue, 0 }  ,fill opacity=1 ][line width=1.5]    (313,204.4) -- (280.49,228) ;
%Straight Lines [id:da9620132957902131]
\draw [fill={rgb, 255:red, 0; green, 0; blue, 0 }  ,fill opacity=1 ][line width=1.5]    (280.49,228) -- (230,228) ;

% Coordinate system
%Straight Lines [id:da6545402950218887]
\draw    (379,257) -- (429,256.04) ;
\draw [shift={(431,256)}, rotate = 538.9] [color={rgb, 255:red, 0; green, 0; blue, 0 }  ][line width=0.75]    (10.93,-3.29) .. controls (6.95,-1.4) and (3.31,-0.3) .. (0,0) .. controls (3.31,0.3) and (6.95,1.4) .. (10.93,3.29)   ;
%Straight Lines [id:da8236119241584237]
\draw    (379,257) -- (379,229) ;
\draw [shift={(379,227)}, rotate = 450] [color={rgb, 255:red, 0; green, 0; blue, 0 }  ][line width=0.75]    (10.93,-3.29) .. controls (6.95,-1.4) and (3.31,-0.3) .. (0,0) .. controls (3.31,0.3) and (6.95,1.4) .. (10.93,3.29)   ;
%Straight Lines [id:da19573210473559044]
\draw    (379,257) -- (411.39,233.18) ;
\draw [shift={(413,232)}, rotate = 503.67] [color={rgb, 255:red, 0; green, 0; blue, 0 }  ][line width=0.75]    (10.93,-3.29) .. controls (6.95,-1.4) and (3.31,-0.3) .. (0,0) .. controls (3.31,0.3) and (6.95,1.4) .. (10.93,3.29)   ;

% Masshilfslinien
%Straight Lines [id:da029740238179323586]
\draw    (179,271) -- (179,235) ;
%Straight Lines [id:da7650347500975532]
\draw    (320.3,271) -- (320.3,235) ;
%Straight Lines [id:da018167774324994168]
\draw    (153,235) -- (179,235) ;
%Straight Lines [id:da37919086237682864]
\draw [line width=0.75]    (153,194) -- (179,194) ;
%Straight Lines [id:da6552493400878835]
\draw [line width=0.75]    (247.2,127.05) -- (273.2,127.05) ;
%Straight Lines [id:da8826068747992396]
\draw [line width=0.75]    (328,203) -- (293,208) ;
%Straight Lines [id:da6085265387648164]
%Arrowed line 1
\draw [line width=0.75]    (317.3,253) -- (182,253) ;
\draw [shift={(179,253)}, rotate = 360] [fill={rgb, 255:red, 0; green, 0; blue, 0 }  ][line width=0.08]  [draw opacity=0] (8.93,-2.29) -- (0,0) -- (8.93,2.29) -- cycle    ;
\draw [shift={(320.3,253)}, rotate = 180] [fill={rgb, 255:red, 0; green, 0; blue, 0 }  ][line width=0.08]  [draw opacity=0] (8.93,-2.29) -- (0,0) -- (8.93,2.29) -- cycle    ;
%Straight Lines [id:da5712161796174655]
\draw    (159,197) -- (159,232) ;
\draw [shift={(159,235)}, rotate = 270] [fill={rgb, 255:red, 0; green, 0; blue, 0 }  ][line width=0.08]  [draw opacity=0] (8.93,-2.29) -- (0,0) -- (8.93,2.29) -- cycle    ;
\draw [shift={(159,194)}, rotate = 90] [fill={rgb, 255:red, 0; green, 0; blue, 0 }  ][line width=0.08]  [draw opacity=0] (8.93,-2.29) -- (0,0) -- (8.93,2.29) -- cycle    ;
%Shape: Rectangle [id:dp8685246838859403]
\draw  [pattern=_ax055x350,pattern size=6pt,pattern thickness=0.75pt,pattern radius=0.75pt, pattern color={rgb, 255:red, 155; green, 155; blue, 155}][dash pattern={on 5.63pt off 4.5pt}][line width=0.5]  (223,201) -- (273.51,201) -- (273.51,231) -- (223,231) -- cycle ;
%Straight Lines [id:da42390508024186657]
\draw [line width=1.5]    (348,215) -- (320.3,235) ;
%Straight Lines [id:da47822940043559947]
\draw [line width=1.5]    (414.5,168.05) -- (373,198) ;

% Text Node
\draw (435,249.4) node [anchor=north west][inner sep=0.75pt]    {$x$};
% Text Node
\draw (417,222.4) node [anchor=north west][inner sep=0.75pt]    {$y$};
% Text Node
\draw (385,216.4) node [anchor=north west][inner sep=0.75pt]    {$z$};
% Text Node
\draw (224,258.4) node [anchor=north west][inner sep=0.75pt]    {$\SI{0.0315}{\meter}$};
% Text Node
\draw (106,206.4) node [anchor=north west][inner sep=0.75pt]    {$\SI{0.0075}{\meter}$};
% Text Node
\draw (148,145.4) node [anchor=north west][inner sep=0.75pt]    {$\SI{0.02405}{\meter}$};
% Text Node
\draw (335,195) node [anchor=north west][inner sep=0.75pt]    {$\Gamma _{spline}$};
% Text Node
\draw (303,167.4) node [anchor=north west][inner sep=0.75pt]    {$\Gamma _{slip}$};
% Text Node
\draw (223,201.4) node [anchor=north west][inner sep=0.75pt]    {$\Gamma_{s}$};
% Text Node
\draw (183,196.4) node [anchor=north west][inner sep=0.75pt]    {$\Gamma _{in}$};
% Text Node
\draw (278.2,130.45) node [anchor=north west][inner sep=0.75pt]    {$\Gamma _{out}$};

\end{tikzpicture}
  \caption{
    Simulation domain for single mixing element as a part of the unwound mixing section with $\SI{7.5}{\milli\meter}$ channel depth.
  }
  \label{fig:simDom}
\end{figure}
%We claim that both, the planar approximation of the fluid domain, as well as the chosen dimensions pose admissible approximations since we do not envision our approach as any sort of submodelling.
%Instead, we argue that through significant reorientation mixing will be enhanced regardless of the actual inflow condition.

The boundary conditions are discussed in the following.
From the bore diameter $d=\SI{60}{\milli \meter}$ and a rotational speed of \SI{8.8}{\per \minute}, we obtain a velocity at the top of $v=\SI{0.027646}{\meter \per \second}$.
Prescribing this rotational movement at the top -- known as barrel rotation or inverse kinematics -- might appear intuitively questionable as the resulting flow field does not seem to resemble physical conditions.
However, we follow the argumentation of Habla et al.\ \cite{Habla13} who prove frame invariance of the melt temperature and in general discuss the negligibility of centrifugal and Coriolis forces.
%Considering further the negligibility of the flow component in \(x\)-direction, we do not cover the mass flow explicitly.
At the inflow, we thus prescribe a linear velocity profile \(\mathbf{u}=\left(0,\hat{v},0\right)\), with \(\hat{v} = 0.027646\frac{z}{0.0075} \si{\meter \per \second}\).
We apply slip boundary conditions on both sides of the flow channel, no-slip on the bottom and mixing element, and a stress-free boundary condition at the outflow.
For the temperature, we apply Dirichlet boundary conditions on both screw (bottom) and barrel (top) with \SI{494}{\kelvin}, and \SI{488}{\kelvin} respectively.

The material properties correspond to $T^*=\SI{473}{\kelvin}$, and a reference temperature of $T_s = \SI{237}{\kelvin}$ is specified.
The material is characterized by density $\rho=\SI{736}{\kilogram \per \cubic \meter}$, specific heat capacity $c_p = \SI{2900}{\joule \per \kilogram \per \kelvin}$, and thermal conductivity $\lambda=\SI{0.4}{\watt \per \meter \per \kelvin}$.
Finally, we prescribe $A=\SI{9,472}{\pascal \second}$, $B=\SI{0.1871}{\second}$, and $C=0.655$ as the Carreau model paramters.

%\todo[inline]{Check that this subsection can really go and if so delete}
%\input{sec_examples_and_validation}
\section{Numerical examples}
\label{sec:results}
In this section, we present a validation of the objective function followed by a verification of the complete design approach.
We start by assessing the objective function's monotonicity and sensitivity -- two main requirements for successful optimization -- before we discuss mesh convergence.
We then show numerical results and finally compare them to experimental data.

\subsection{Objective evaluation}
As stated in Sec.~\ref{subsec:goveq_mixing}, we aim to quantify mixing based on the interfacial area.
The interfacial area is computed as the summed-up lengths of the convex hulls around  $m \times n$ particle sets.
To define these particles sets, we first segment the inflow into  $m \times n$ rectangular subdomains.
Then, we track the particles covered by each of these subdomains and construct their convex hulls again at the outflow.
Comparing each convex hull's length between inflow and outflow, finally, gives a measure of mixing.
Assuming $2 \times 2$ rectangular inflow subdomains, possible outflow convex hulls are shown in Fig.~\ref{fig:three_objectives}.
Convex hulls, however, might be only an approximation of the actual particle sets' boundaries.
To understand its influence, we start by discussing two primitive flow patterns, and their influence on the initially rectangular convex hulls.\par
First, we consider a flow skew to the computational domain.
Such a flow will simply shift all particles uniformly and will thus not cause mixing.
Similarly, all convex hulls are uniformly displaced and thus -- with no distortion between in- and outflow -- no indication of increased mixing is expected.
%It should be noted that if particle's displacements were computed using trivial distance measures like the $L_1$ or $L_2$ norm -- as opposed to the convex hull -- these measures fail in this case and indicate enhanced mixedness of the flow.
A second simple scenario is shear, i.e., Couette flow.
In Couette flows, the linear velocity profile leads to linearly varying particle displacements and thus increased mixing.
Again, the convex hull will follow the particle displacements.
However, now, it will be deformed into a parallelogram shape.
The parallelogram features an increased perimeter length, and thus correctly indicates increased mixing.\par
%These simple cases indicate the applicability of our objective function in contrast trivial distance measures like the $L_1$ or $L_2$ norm applied to each particle.
These two scenarios indicate the general applicability of our method for cases in which the particle sets' boundaries are indeed convex.
To further prove the objective's suitability, we extend our analysis to non-convex situations.
We obtain non-convex particle sets using an artificial example in which we mimic the influence of two counter-rotating vortices in the mixing element's wake.
%Obtained from these tests, Fig.~\ref{fig:three_objectives} shows examples of convex hulls at the outflow domain, which are obtained by introducing two counter-rotating vortices.
These vortices are modeled as angular displacements.
To resemble physically reasonable displacements, a smooth transition to the far field is ensured by scaling the imposed displacement as $\phi=\Delta \phi \max\left(0,\cos(r)\right)$ as illustrated in Fig.~\ref{fig:DisplacementExplanation}.
\begin{figure}[!h]
	\centering
	\input{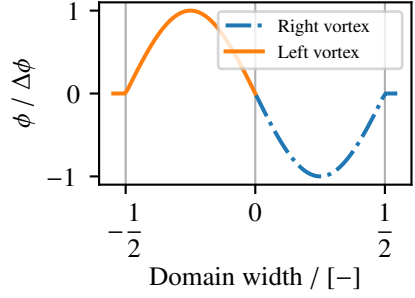}
	\caption{Schematic view of the imposed angular displacement imposed by the two vortices plotted over the horizontal domain symmetry line.}
	\label{fig:DisplacementExplanation}
\end{figure}
In particular, $r$ is chosen as $\min\left(\frac{1}{2}h, \frac{1}{2} w\right)$, and $h$ and $w$ respectively denote the height and width of the shown section.
A study with three angular offsets of varying degrees is shown in Fig.~\ref{fig:three_objectives}.
In view of the increasing objective function value for increasing displacemens, it is clearly visible that the approximate nature of the convex hull appears negligible for small offsets. 
\begin{figure}[!h]
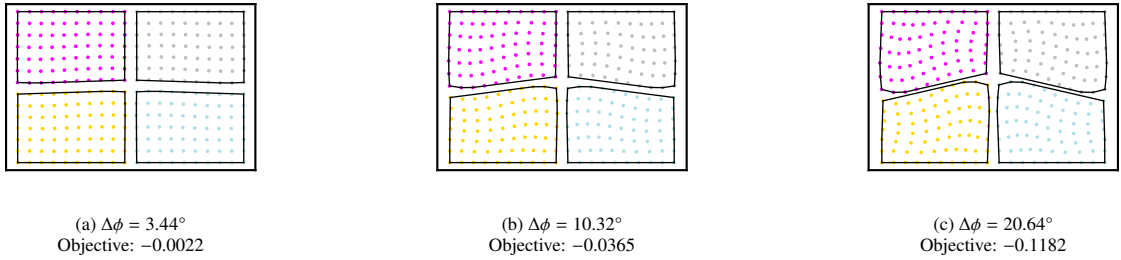

     \centering
     \begin{subfigure}[b]{0.31\textwidth}
         \centering
         \input{figures/fig_5/000.png.pgf}
             \caption{$\Delta \phi = \SI{3.44}{\degree}$ \\ Objective: $-0.0022$}
         \label{fig:obj_sample1}
     \end{subfigure}
     \hfill
     \begin{subfigure}[b]{0.31\textwidth}
         \centering
         \input{figures/fig_5/003.png.pgf}
             \caption{$\Delta \phi = \SI{10.32}{\degree}$ \\ Objective: $-0.0365$}
         \label{fig:obj_sample2}
     \end{subfigure}
     \hfill
     \begin{subfigure}[b]{0.31\textwidth}
         \centering
         \input{figures/fig_5/006.png.pgf}
             \caption{$\Delta \phi = \SI{20.64}{\degree}$ \\ Objective: $-0.1182$}
         \label{fig:obj_sample3}
     \end{subfigure}
        \caption{Sample convex hulls at $\Gamma_{out}$ constructed for a 2x2 split of the inflow domain $\Gamma_s$.
        Artificial displacements are introduced as rotational offsets, i.e., we first transform all points into polar coordinates and then impose the additional offset $\phi$.
        For the left half of the domain, we construct polar coordinates around $\left(\frac{1}{4}w,\frac{1}{2}h\right)$, whereas for the right half we choose $\left(\frac{3}{4}w,\frac{1}{2}h\right)$ which denote the corresponding rotation centers.
     %   To ensure a smooth transition towards the boundary, the actual displacement is decayed by $\phi=\Delta \phi \max\left(0,\cos(r)\right)$, where $r$ is chosen as $\min\left(\frac{1}{2}h, \frac{1}{2} w\right)$, with $h$ and $w$ being the height and width of the shown outflow region respectively.
        }
        \label{fig:three_objectives}
\end{figure}

%For a visualization, all cases of planar shear might be viewed as a one-dimensional simplification of the displacements depicted in %Fig.~\ref{fig:three_objectives} and are thus explicitly not shown here.
%It should be noted though, that simple shear is correctly correlated to increased mixing by the proposed objective.
%The same argumentation is applied considering multiple or a single rotational center, indicating that the proposed objective is robust %to simple redirections yet sensitive to all forms of shear.
After studying the objective function for convex and slightly non-convex particle distributions, the remaining question is how well the objective function behaves for more pronounced flow reorientations.
Figure \ref{fig:objTest} illustrates the results obtained from the vortex example for rotation angles $\Delta \phi$ up to $\SI{103.2}{\degree}$.
\begin{figure}[!h]
  \centering
  \input{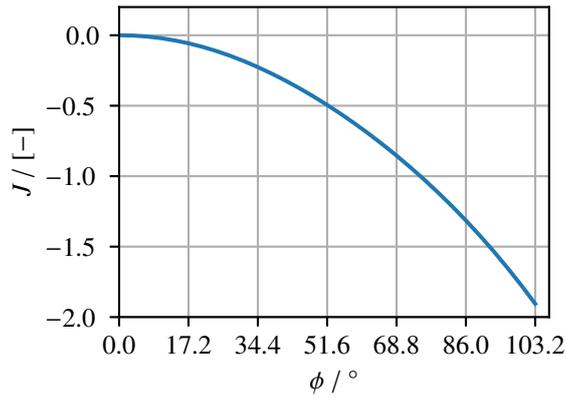}
  \caption{Monotonicity and sensitivity of chosen objective shown for two counterrotating centers as exemplified in Fig.~\ref{fig:three_objectives}.}
  \label{fig:objTest}
\end{figure}
It is seen that even for such strong vortices, the proposed objective function displays sensitivity and monotonic behavior.\par
Further tests using chaotic, i.e., random, flows are not conducted in view of the flow's laminarity.
Thus, we claim that the proposed objective function is not only justifiable -- built on experimental evidence -- but also numerically advantageous.

\subsection{Mesh convergence study}
To define a suitable computational mesh, we study mesh convergence by applying a sample configuration with $\alpha=10\degree$, $\beta=15\degree$, $\gamma=45\degree$, height scaling of $h=0.9$, a shrink ratio of $s=1.9$ and zero global rotation around $z$.
We show convergence plots for the solution of our governing equations, but also of the objective function.
More specifically, we focus on the objective function as the actual measure driving our optimization pipeline.\par
Since no analytical solution is available, we construct a very fine mesh ($2.1\times 10^6$ elements) and compare the nodal relative root mean square error of the projected solution against this reference.
To obtain a constant number of particle traces, the resolution at $\Gamma_{s}$ and $\Gamma_{out}$ is kept constant.
\begin{figure}[!h]
  \centering
  \input{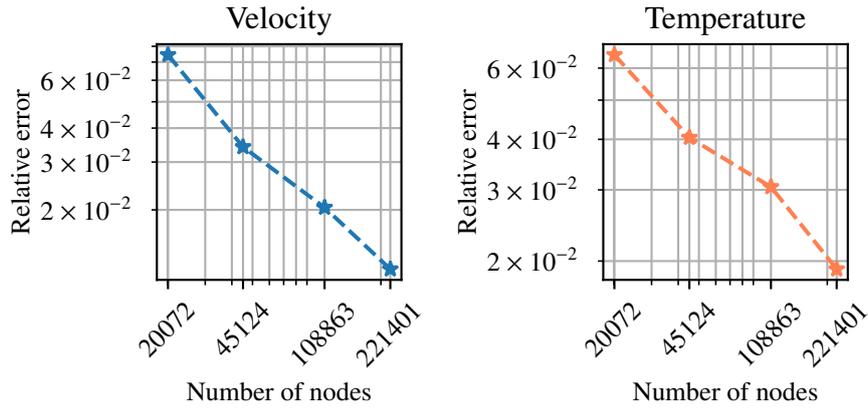}
  \caption{
    Relative error of velocity and temperature compared to the finest mesh.
  }
  \label{fig:meshStudy}
\end{figure}
As can be seen in Fig.~\ref{fig:meshStudy}, even for the coarse mesh, the relative error in the flow variables is below $\SI{8}{\percent}$.
Contrary to the small errors for velocity and temperature, the objective function shows a significant dependency on the mesh resolution as depicted in Fig.~\ref{fig:meshStudyObj}.
In particular, the coarsest mesh seems unsuited.
\begin{figure}[!h]
	\centering
	\input{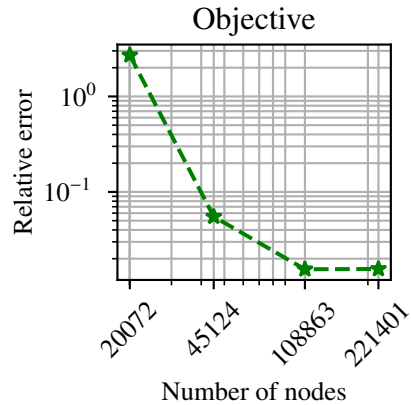}
	\caption{
		Relative error of the objective value compared to the finest mesh.
	}
	\label{fig:meshStudyObj}
\end{figure}
However, despite mesh convergence for velocity and temperature, the objective error levels out with increased mesh resolution.
The observed plateau indicates that an even finer resolution of $\Gamma_{out}$ or higher-order interpolation would be needed to track the particles more precisely.
However, for the numerical studies, we assume that the obtained accuracy suffices and choose the mesh with 108,863 nodes and 586,076 elements.

\subsection{Optimization and experimental verification}
\label{subsec:exp_res}
In this section we present the verification of the proposed design framework.
To verify our approach, we first shape optimize one example mixing element.
Then, an existing mixing section is combined with the obtained optimized mixing element, and finally mixing is assessed.\par
The chosen mixing element's initial geometry is shown in the top left corner of Fig.~\ref{fig:geos}.
% 
%simulation setup, the splitting of the optimization problem, and the utilized objective function, the optimization process is applied to a standard rhomboid mixing element shown in the top left corner of Fig.~\ref{fig:geos} using the shape
%optimization library Dakota \cite{dakota}.
%
The shape optimization is driven by the Dakota library \cite{dakota}.
Dakota provides a variety of optimization algorithms, of which the ncsu\_direct algorithm is chosen for its global convergence.
After 193 iterations, an optimal configuration is found with an objective value of $J=-0.1078$.
It should be noted that this objective function value has no direct physical interpretation.
This optimized mixing element, as well as the assembled mixing sections used in the experimental verification, are also shown in Fig.~\ref{fig:geos}.
\begin{figure}[!h]
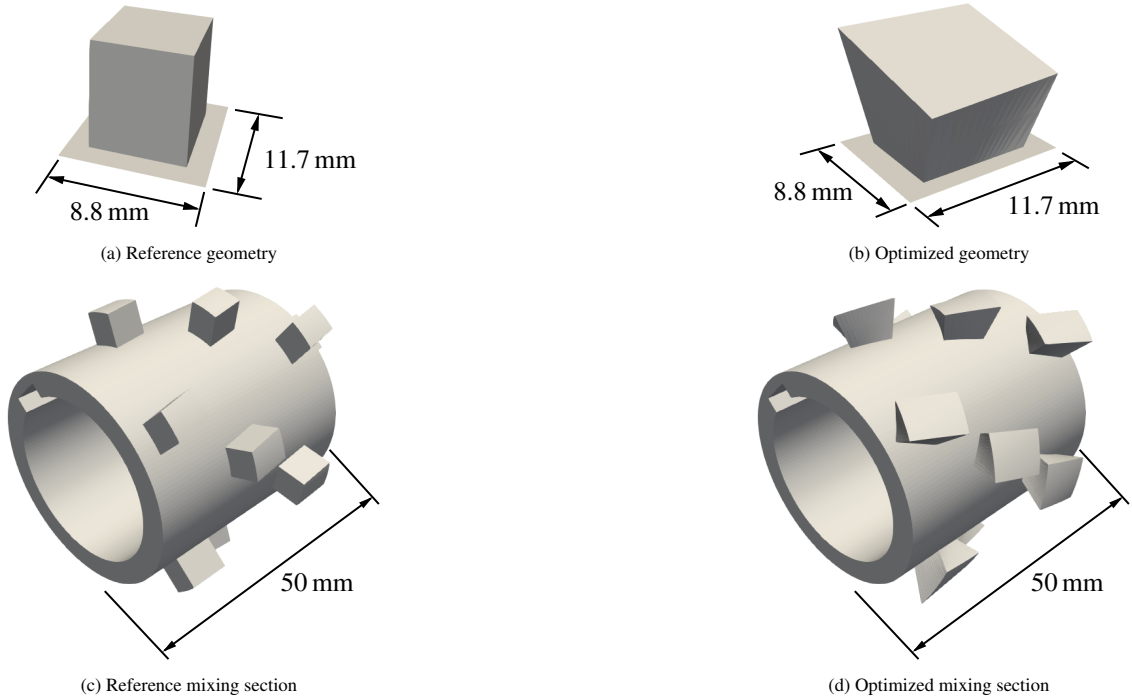

     \centering
     \begin{subfigure}[t]{0.4\textwidth}
         \centering
             \input{figures/fig_9/elem_ref.tex}
         \caption{Reference geometry}
         \label{fig:ref_geo}
     \end{subfigure}
     \hfill
     \begin{subfigure}[t]{0.4\textwidth}
         \centering
             \input{figures/fig_9/elem_opt.tex}
         \caption{Optimized geometry}
         \label{fig:opt_geo}
     \end{subfigure}
     \\
     \begin{subfigure}[b]{0.4\textwidth}
         \centering
             \input{figures/fig_9/screw_ref.tex}
         \caption{Reference mixing section}
         \label{fig:ref_geo}
     \end{subfigure}
     \hfill
     \begin{subfigure}[b]{0.4\textwidth}
         \centering
         \input{figures/fig_9/screw_opt.tex}
         \caption{Optimized mixing section}
         \label{fig:opt_geo}
     \end{subfigure}
        \caption{Comparison of initial and optimized geometry.
        }
        \label{fig:geos}
\end{figure}

Using the the initial and the numerically designed mixing sections, we verify our approach experimentally.
Therefore, we measure the additive and temperature distribution using the experiment described by Hopmann et al. \cite{hopmann2020method}.
The additive distribution is evaluated at the nozzle using the average grey-scale value, $\mu_{grey}$ .
A lower value of  $\mu_{grey}$  indicates increased homogeneity.\par
Given the laminarity of the flow and negligible effects of diffusion and heat conduction \cite{Mohr57}, it is expected that advection is also the dominant mechanism for heat exchange and thus -- despite not being explicitly formulated in the objective -- 
improvements towards a homogeneous temperature distribution will result from our approach.
The experimental data, as shown in Figures~\ref{fig:expRes} and \ref{fig:expTempRes}, confirms this assumption and shows good agreement with numerical results.
Up to \SI{5}{\degreeCelsius} temperature reduction and \SI{25}{\percent} increase in distributive mixing is measured using the shape-optimized mixing element.
It should be noted that the experimental results  reveal surprisingly satisfying benefits in the off-design region as well.
\begin{figure}[!h]
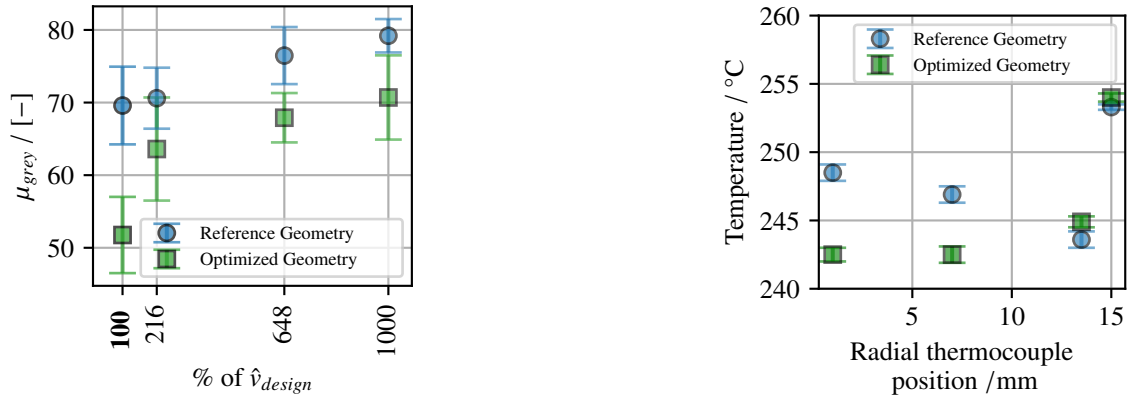

	\centering
	\begin{subfigure}[t]{0.49\textwidth}
		\input{figures/fig_10/expResults.pgf}
		\caption{Distributive mixing determined from experimental studies.
			The optimized geometry performs better not only in the desired, but also for off-design mass flows indicating the ability of the proposed simple objective to isolate and emphasize physically meaningful shape improvements.
		}
		\label{fig:expRes}
	\end{subfigure}
	\hfill
	\begin{subfigure}[t]{0.49\textwidth}
		\centering
		\input{figures/fig_10/expTempResults.pgf}
		\caption{Temperature profiles at different radii determined with a thermocouple at \SI{1000}{\percent} $\hat{v}_{design}$.
			Consistent with \ref{fig:expRes}, the optimized geometry provides enhanced homogenization of the temperature despite the fixed barrel temperature.
		}
		\label{fig:expTempRes}
	\end{subfigure}
	\caption{Experimental results obtained from reference and optimized geometry.}
\end{figure}

We study temperature distributions at higher speeds than the design point.
This is a consequence of the low design speed which causes high residence times, and thus the temperature distributions are strongly influenced by the extruder's heat regulation system.\par
The usage of  $\mu_{grey}$ (the average of all grey-scale values) in place of the typical $\delta_{grey}$ (the standard deviation of the grey-scale values) stems from the fact that the experimentally determined values of $\delta_{grey}$ do not reflect the visually apparent differences between the samples.
Figure \ref{fig: grey_scale_values} provides a more detailed view: The blue curves depict the distribution of grey-scale values for seven extrudate samples processed with the reference geometry.
In contrast, the green curves show the grey-scale values for  seven extrudate samples processed with the optimized mixer. A grey-scale value of $0$ shows a purely black pixel, a value of $255$ a purely white pixel.
Above the plots, a representative extrudate sample is shown for both geometries. 
\begin{figure}[!h]
	\centering
	\input{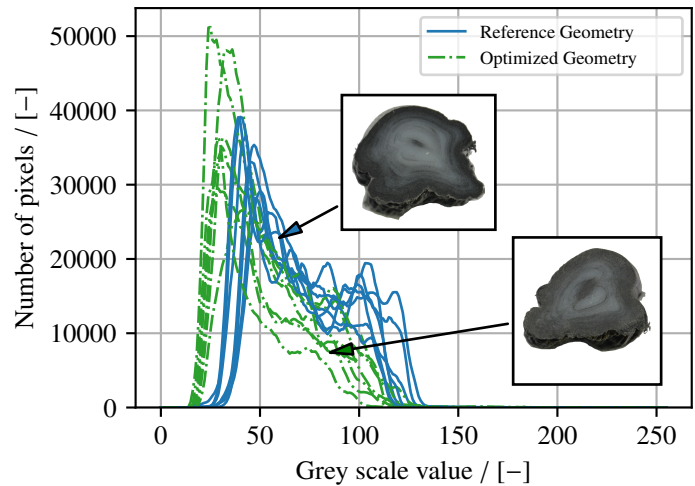}
	\caption{Comparison of grey-scale distributions for reference and optimized geometry at the design point for seven different samples with exemplary cross section images.
	For a detailed description of the experimental setup, the reader is referred to \cite{hopmann2020method}.}
	\label{fig: grey_scale_values}
\end{figure}
It is apparent that samples processed with the optimized geometry show a wider black brim along the circumference. In addition, radially distanced concentric black stripes appear on the samples processed with the optimized geometry, but not on the reference samples. 
This is the consequence of the increase in radial flow in the optimized geometry.
%\begin{figure}[!h]
%	\centering
	%\includegraphics[width=1.0\textwidth]{Media/BILD_IKV_SEITE.png}
	%\input{IKV_Validation/GreyValuesFromPhotos/greyscaleTest.pgf}
%	\input{greyscaleTest.pgf}
%	\caption{Comparison of grey-scale distributions for reference and optimized geometry at the design point for seven different samples with exemplary cross section images.}
%	\label{fig: grey_scale_values}
%\end{figure}
However, the grey-scale distribution fails to become broader with the introduction of the optimized geometry. 
Instead, it appears shifted.
Consequently, $\delta_{grey}$ is not expected to change, while $\mu_{grey}$  drops. 
In practical terms, this represents darker extrudate samples.
As the amount of color batch added remained constant during the experiments, the darkening represents a more efficient spreading of the color batch which is very much in agreement with the definition of good distributive mixing.
\section{Conclusions and outlook}
\label{sec:conclusion}
In the presented work, a finite-element-based shape optimization framework is developed.
Using PDE-based mesh update methods, FFD, and a novel, experiment-inspired objective function based on interfacial area measurements, the resulting framework is used to optimize existing mixing elements.
In addition, a verification of the proposed approach against experimental evidence is presented.\par
As in other disciplines, numerical design methods, again, proved to be a valuable design aid capable of identifying also counterintuitive shapes.
In fact, the chosen framework repeatedly produces mixing elements with inclined top sections, as shown in Fig.~\ref{fig:geos}.
These reproducible results are unexpected in the sense that such mixing elements are not currently used in industrial practice.
However, review of the very early literature on laminar mixing \cite{Erwin78_TheoryOfLaminarMixing} reveals that the new design is likely to produce flow patterns similar to\textit{ pure} instead of \textit{simple shear} flow.
Simple shear flow -- obtained in, e.g., Couette flows -- is significantly inferior in enhancing mixing compared to pure shear -- which, being connected to elongation, occurs in tapering flow sections.
With its inclined top, however, the shape-optimized mixing element features such a tapered design.\par
Given the overall good agreement of experimental and numerical results, and the unintuitive solutions obtained, we claim that the benefits of numerical design are demonstrated.
Now, the logical extension towards an overall enhanced design of single-screw-extruder screws is to place a suitable number of mixing elements relative to each other such that the mixing becomes even more pronounced.
We expect that a gradual phase shift in the sense of rotational offsets will be most efficient in providing better mixing.
It should also be noted that the framework's applicability towards the optimization of dispersive mixing, i.e., the reduction of agglomerates, is possible without modification, given a suitable objective function.\par
Finally, we would like to emphasize that the presented finding of tilted-top mixing elements might be crucial for industrial applications since -- assuming the tilted top is the mixing booster -- manufacturability is straightforward, but to the authors' knowledge, such a geometry is not utilized yet.

\section{Acknowledgements}
We gratefully acknowledge the funding received from the German Research Foundation (DFG) via the DFG grant “Automated design and optimization of dynamic mixing and shear elements for single-screw extruders” (HO 4776/37-1 and EL 745/5-1).
Implementation and computations were carried out on the high-performance computing resources provided by IT Center at RWTH Aachen University.

%\section*{References}

\bibliography{mybibfile}

\end{document}